\documentclass{amsart}
\usepackage{amsmath}%
\usepackage{amsfonts}%
\usepackage{amssymb}%
\usepackage{graphicx}
\theoremstyle{plain}

\numberwithin{equation}{section}
\begin{document}
\title[Probabilistic Weapon-Target Allocation Models]{Probabilistic Weapon-Target Allocation Models for the Critical Infrastructure Protection}
\author{Vladimir NORKIN}
\address[V. Norkin]
{V.M. Glushkov Institute of Cybernetics, Kyiv, Ukraine \newline%
\indent National Academy of Sciences of Ukraine}%
\email[V. Norkin]{vladimir.norkin@gmail.com}%
\urladdr{https://incyb.kiev.ua/en/employee/norkin-volodimir-ivanovich}
\date{November 07, 2026}

\begin{abstract}
This paper analyzes methods for solving optimal resource/weapon allocation problems for the protection of critical infrastructure. It examines nonlinear discrete "Attacker-Defender" and "Weapon-Target Allocation" models, incorporating criteria based on the difference in expected losses and risk-related constraints. Methods are proposed for solving these problems by reformulating them as mixed-integer linear programming problems.
\end{abstract}
\maketitle

\section{Introduction}

The Weapon-Target Allocation (WTA) problem is a classic problem in operations research \cite{manne1958target}, representing a generalization of the optimal assignment problem (assigning agents to tasks). A detailed overview of the Weapon-Target Assignment (WTA) problem is provided in \cite{li2024comprehensive,liu2025current}, while dynamic models are discussed in \cite{cai2006survey,hasan2020weapon,li2023weapon,okumura2023solving,xin2010efficient}. Early publications concerning the Missile-Allocation Problem (MAP, see \cite{ahuja2007exact}) are analyzed in
\cite{matlin1970review,gurin1968problems}, which also offer a classification of such problems based on their characteristic features. Reference \cite{lu2021new} presents a reformulation of the WTA problem as a large-scale problem involving binary variables; the authors propose a fast, exact algorithm for its solution and provide comprehensive results from computational experiments. Work \cite{bertsimas2025solving} investigates very large-scale problems, considers robust and stochastic formulations, and accounts for uncertainty in probabilistic parameters. Reference \cite{kline2019weapon} provides an overview of nine static and ten dynamic models (including two-stage, multi-stage, and game-theoretic formulations) for the WTA problem, and describes various exact, heuristic, and metaheuristic algorithms for solving both static and dynamic WTA problems. Finally, \cite{andersen2022weapon} presents a literature review and proposes an approach for reducing very large-scale problems to mixed-integer linear programming problems using linearization techniques.

This paper examines new spatial weapon-target allocation problems within a probabilistic framework, where both weapons and targets possess spatial coordinates (an additional spatial index). The optimization criteria include not only the value of the engaged targets but also the probability of their successful engagement. This approach can also be used to model the defense of one or more critical infrastructure assets against aerial attacks. In this context, the model accounts for the value of intercepted aerial targets, the cost of the weapons employed, and the value of any damage sustained by the assets. Similarly, problem formulations aimed at minimizing the probability of damage to critical infrastructure assets are considered. These formulations represent large-scale, non-linear integer optimization problems. To solve them, the literature proposes the use of both exact and heuristic methods.

In this work, the problems are first equivalently transformed—by introducing auxiliary variables—into mixed-integer programming problems with convex continuous objective functions and linear constraints. Subsequently, the objective functions are approximated from above and below, with arbitrary precision, using convex functions defined as the maximum of a set of linear functions. Similar piecewise-linear approximations, in various forms, have been employed in numerous studies (see \cite{manne1958target,ahuja2007exact,andersen2022weapon,lotter2013multiobjective,camm2002nature}). The resulting problems—minimizing a maximum function subject to linear constraints—can be straightforwardly reduced to mixed-integer linear programming problems, with the precision of the approximations being controllable. Furthermore, additional constraints are derived for all auxiliary continuous variables, which serves to restrict the feasible region of the problems. The resulting approximate problems are solved using standard mixed-integer linear programming methods.

Furthermore, the standard formulations of the Weapon-Target Assignment (WTA) problem are extended to include non-linear constraints on the risk (probability) of individual targets remaining unhit; non-cost-based formulations, in which the objective function represents the probability of targets remaining unhit, are also considered. However, the proposed approach for linearizing these problems remains applicable to such modifications as well.

\section{Nomenclature of WTA problems}

Mathematically, this WTA problem is formulated as follows \cite{manne1958target}:
\begin{equation}
\sum\nolimits_{i\in T}{a}_{i}\left( 1-\prod\nolimits_{j\in W}\left( 1-{p}_{ij} \right)^{{{x}_{ij}}} \right)\to {\max }_{\left\{ {x}_{ij} \right\}},\;\;\; \label{(1)} 
\end{equation}
(maximization of expected damage),
\begin{equation}
\sum\nolimits_{i\in T}{{{x}_{ij}}}\le {{D}_{j}}, \;\;\;j\in W,\;\;\;  \label{(2)}
\end{equation}
(restriction on the availability of weapons),
\begin{equation}
{x}_{ij}\in {\mathbb{Z}_{+}}, \;\;\;i\in T,\ j\in W, \;\;\; \mbox{(restriction on integers)}, \label{(3)}
\end{equation}
where $T$ denotes the set of targets, among which there may be targets of the same type; $W$- the set of weapons; ${{a}_{i}}$ - the cost of target $i\in T$; ${{p}_{ij}}$ - the probability of hitting target $i$ with weapon $j$; ${{x}_{ij}}$ - the number of weapons of type $j$ that are simultaneously (in a volley) used against target $i$; ${{D}_{j}}$ - the available number of weapons of type $j$; ${{\mathbb{Z}}_{+}}$ - the set of non-negative integers. Here \[{{\pi }_{i}}=\prod\nolimits_{j\in W}{{{\left( 1-{{p}_{ij}} \right)}^{{{x}_{ij}}}}}=\prod\nolimits_{j\in W}{{{\left( {{q}_{ij}} \right)}^{{{x}_{ij}}}}}\] is the probability of not hitting target $i$ with simultaneous independent use (volley) of ${{x}_{ij}}$ units of weapon type $j$ against target $i$; ${{q}_{ij}}=1-{{p}_{ij}}$. If for some $i$ all ${{x}_{ij}}=0,\,\,j\in W,$ then this target is not attacked. We will denote $\left| T \right|$ and $\left| W \right|$ the number of elements in the sets $T$ and $W$ . Problem \eqref{(1)}--\eqref{(3)} is a nonlinear integer programming problem.

Weapons are lost after use and can be expensive. Let ${{b}_{j}}\ge 0$ be the cost of firing a weapon of type $j\in W$. Then the cost of the weapon used is $\sum\nolimits_{i\in T,j\in W}{{{b}_{j}}{{x}_{ij}}}$. This is an additional objective function for problem \eqref{(1)}--\eqref{(3)}. Thus, the problem becomes two-criteria and the question arises of constructing the Pareto-optimal frontier and non-dominated solutions of the two-criteria problem.

Other (binary) formulations of WTA problems \cite{manne1958target,lu2021new,kline2019weapon} have the following form:

minimize $\sum\limits_{i\in T}{{{a}_{i}}\prod\nolimits_{j\in W}{\left( 1-{{p}_{ij}}{{x}_{ij}} \right)}}$

under the constraints: $\sum\limits_{i\in T}{{{x}_{ij}}}\le 1$, $j\in W$;

${{x}_{ij}}\in \left\{ 0,1 \right\}$, $i\in T$, $j\in W$.

Another binary objective function in this problem has the form \cite{ahuja2007exact}:

	$\sum\limits_{i\in T}{{{a}_{i}}\prod\nolimits_{j\in W}{{{\left( 1-{{p}_{ij}} \right)}^{{{x}_{ij}}}}}}=\sum\limits_{i\in T}{{{a}_{i}}\prod\nolimits_{j\in W:\,{{x}_{ij}}=1}{\left( 1-{{p}_{ij}} \right)}}$.
	
	In the WTA problem taking into account the cost of weapons, the objective function has the form \cite{ohanley2013probability}:
	
	$\sum\limits_{i\in T}{{{a}_{i}}\prod\nolimits_{j\in W}{{{\left( 1-{{p}_{ij}} \right)}^{{{x}_{ij}}}}}}+\sum\limits_{i\in T,j\in W}{{{b}_{j}}{{x}_{ij}}}$.
		
To construct the Pareto-optimal frontier, it is necessary to solve a series of problems with a variable parameter $B\in \left[ 0,\sum\nolimits_{j\in W}{{{b}_{j}}\left| {{D}_{j}} \right|} \right]$ of permissible losses in the cost of the applied weapons. To do this, constraints of the form are added to the constraints of problems \eqref{(1)}--\eqref{(3)} or \eqref{(11)}--\eqref{(14)}:
	\[\sum\nolimits_{i\in T,j\in W}{{{b}_{j}}{{x}_{ij}}}\le B.\]
	
	Sometimes, instead of a two-criteria formulation, the problems of optimal distribution of weapons among targets are considered using the criterion of the difference of losses \cite{perevozchikov2019problem}:
\begin{equation}\label{(4)}
\left[ \sum\nolimits_{i\in T}{{{a}_{i}}}\left( 1-{{\prod\nolimits_{j\in W}{\left( 1-{{p}_{ij}} \right)}}^{{{x}_{ij}}}} \right)-\sum\nolimits_{i\in T,j\in W}{{{b}_{j}}{{x}_{ij}}} \right]\to {{\max }_{\left\{ {{x}_{ij}} \right\}}},
\end{equation}	
\begin{equation}\sum\nolimits_{i\in T}{{{x}_{ij}}}\le {{D}_{j}},	\;\;\;j\in W,\label{(5)}\end{equation}	
\begin{equation}{{x}_{ij}}\in {{\mathbb{Z}}_{+}},\;\;\;	i\in T,\ j\in W.\label{(6)}\end{equation}	

Problem \eqref{(4)}--\eqref{(6)} can be equivalently reformulated as a minimization problem:
\begin{equation}\label{(7)}
\left[ \sum\nolimits_{i\in T}{{{a}_{i}}}{{\prod\nolimits_{j\in W}{\left( 1-{{p}_{ij}} \right)}}^{{{x}_{ij}}}}+\sum\nolimits_{i\in T,j\in W}{{{b}_{j}}{{x}_{ij}}} \right]\to {{\min }_{\left\{ {{x}_{ij}} \right\}}}
\end{equation}
\begin{equation}
\sum\nolimits_{i\in T}{{{x}_{ij}}}\le {{D}_{j}},	\;\;\;j\in W, \label{(8)}
\end{equation}
\begin{equation}
{{x}_{ij}}\in {{\mathbb{Z}}_{+}}, \;\;\;i\in T,\ j\in W.\label{(9)}
\end{equation}

Note (about consecutive shots). With consecutive (maximum ${{x}_{ij}}$) shots at target $i$ by means $j$ until hit, the probability of hit is
${{P}_{ij}}={{p}_{ij}}+(1-{{p}_{ij}}){{p}_{ij}}+...+{{(1-{{p}_{ij}})}^{{{x}_{ij}}-1}}{{p}_{ij}}=1-{{(1-{{p}_{ij}})}^{{{x}_{ij}}}}$,
whereas on average, the following number of shots is used:
\begin{eqnarray*}
  {{X}_{ij}}&=&{{p}_{ij}}+2(1-{{p}_{ij}}){{p}_{ij}}+3{{(1-{{p}_{ij}})}^{2}}{{p}_{ij}}+...+{{x}_{ij}}{{(1-{{p}_{ij}})}^{{{x}_{ij}}-1}}{{p}_{ij}} \\ 
 &=&\frac{1}{{{p}_{ij}}}\left( 1-({{x}_{ij}}+1){{(1-{{p}_{ij}})}^{{{x}_{ij}}}}+{{x}_{ij}}{{(1-{{p}_{ij}})}^{{{x}_{ij}}+1}} \right) \\ 
 &=&\frac{1}{{{p}_{ij}}}\left( 1-({{x}_{ij}}{{p}_{ij}}+1){{(1-{{p}_{ij}})}^{{{x}_{ij}}}} \right)<\frac{1}{{{p}_{ij}}}. 
\end{eqnarray*}
If means of destruction $j$ are used sequentially (echeloned defense) in quantities ${{x}_{ij}}$, $j\in W$, then the probability of hitting target $i$ is
\begin{eqnarray*}
  {{P}_{i}}&=&{{P}_{i1}}+(1-{{P}_{i1}}){{P}_{i2}}	+(1-{{P}_{i1}})(1-{{P}_{i2}}){{P}_{i3}}+...\\
	&&+(1-{{P}_{i1}})\times ...\times (1-{{P}_{i(\left| W \right|-1)}}){{P}_{i\left| W \right|}} \\ 
 &=&1-{{(1-{{p}_{i1}})}^{{{x}_{i1}}}}+{{(1-{{p}_{i1}})}^{{{x}_{i1}}}}\left( 1-{{(1-{{p}_{i2}})}^{{{x}_{i2}}}} \right)\\
&&+{{(1-{{p}_{i1}})}^{{{x}_{i1}}}}{{(1-{{p}_{i2}})}^{{{x}_{i2}}}}\left( 1-{{(1-{{p}_{i3}})}^{{{x}_{i3}}}} \right)+... \\ 
 &&+{{(1-{{p}_{i1}})}^{{{x}_{i1}}}}{{(1-{{p}_{i2}})}^{{{x}_{i2}}}}\times ...\times {{(1-{{p}_{i\left( \left| W \right|-1 \right)}})}^{{{x}_{i\left( \left| W \right|-1 \right)}}}}\left( 1-{{(1-{{p}_{i\left| W \right|}})}^{{{x}_{i\left| W \right|}}}} \right),  
\end{eqnarray*}
where $\left| W \right|$ denotes the number of weapons used. The average cost of weapons used is estimated as
\begin{align*}
  & {{b}_{1}}{{X}_{i1}}{{P}_{i1}}+\left( {{b}_{1}}{{x}_{i1}}+{{b}_{2}}{{X}_{i2}} \right)(1-{{P}_{i1}}){{P}_{i2}}+\left( {{b}_{1}}{{x}_{i1}}+{{b}_{2}}{{x}_{i2}}
	+{{b}_{3}}{{X}_{i3}} \right)(1-{{P}_{i1}})(1-{{P}_{i2}}){{P}_{i3}}+... \\ 
 & \;\;\;\;\;\;\;\;\;\;\;\;\;\;\;...+\left( {{b}_{1}}{{x}_{i1}}+{{b}_{2}}{{x}_{i2}}+...+{{b}_{\left( \left| W \right|-1 \right)}}{{x}_{i\left( \left| W \right|-1 \right)}}+{{b}_{\left| W \right|}}{{X}_{i\left| W \right|}} \right)\times\\
&\;\;\;\;\;\;\;\;\;\;\;\;\;\;\;\times(1-{{P}_{i1}})(1-{{P}_{i2}})\times ...
\times (1-{{P}_{i\left( \left| W \right|-1 \right)}}){{P}_{i\left| W \right|}}.  
\end{align*}
The average cost of hit targets alone in \eqref{(1)} may not be sufficient to assess the effectiveness of the weapons distribution, since some targets must be hit with almost certainty. This requirement is taken into account by the additional constraints for problem \eqref{(1)}--\eqref{(3)} that the probability of hitting target $i$ must be at least $1-{{\varepsilon }_{i}}$:

$1-{{\prod\nolimits_{j\in W}{\left( 1-{{p}_{ij}} \right)}}^{{{x}_{ij}}}}\ge 1-{{\varepsilon }_{i}}$
or
${{\prod\nolimits_{j\in W}{\left( 1-{{p}_{ij}} \right)}}^{{{x}_{ij}}}}\le {{\varepsilon }_{i}}$,\;\;\;	
${{\varepsilon }_{i}}\in \left( 0,1 \right]$,  \;\;\; $i\in T$.
If ${{\varepsilon }_{i}}=1$, then the corresponding constraint is satisfied automatically and may not be taken into account in the problem statement. These inequalities can be equivalently rewritten in linear form:
\begin{equation}\label{(10)}
\sum\nolimits_{j\in W}{\left( \ln (1-{{p}_{ij}}) \right)}\,{{x}_{ij}}\le \ln {{\varepsilon }_{i}},\;\;\; 	
i\in T.
\end{equation}

Problem \eqref{(1)}-\eqref{(3)}, \eqref{(10)} can also be equivalently reformulated as a problem of minimizing the cost of unhit targets:
\begin{equation}
\sum\nolimits_{i\in T}{{{a}_{i}}{{\prod\nolimits_{j\in W}{\left( 1-{{p}_{ij}} \right)}}^{{{x}_{ij}}}}}\to {{\min }_{\left\{ {{x}_{ij}} \right\}}}, \label{(11)}
\end{equation}
\begin{equation}
\sum\nolimits_{i\in T}{{{x}_{ij}}}\le {{D}_{j}}, \;\;\;j\in W,	\label{(12)}
\end{equation}
\begin{equation}
	\sum\nolimits_{j\in W}{\left( \ln (1-{{p}_{ij}}) \right)}\,{{x}_{ij}}\le \ln {{\varepsilon }_{i}}, \;\;\;i\in T, \label{(13)}
\end{equation}
\begin{equation}
{{x}_{ij}}\in {{\mathbb{Z}}_{+}}, \;\;\;	i\in T,\ j\in W,  \label{(14)}
\end{equation}
since the constant $\sum\nolimits_{i\in T}{{{a}_{i}}}$ in the objective function \eqref{(1)} does not affect the solution of the problem.

A partial case of problem \eqref{(11)}--\eqref{(14)} is the problem of distributing uniform weapons in quantity $D$ between several targets:
\begin{equation}\label{(15)}
\sum\nolimits_{i\in T}{{{a}_{i}}}{{\left( 1-{{p}_{i}} \right)}^{{{x}_{i}}}}\to {{\min }_{\left\{ {{x}_{i}} \right\}}}
\mbox{\;\;\;(minimization of expected non-injury)},
\end{equation}
\begin{equation}\label{(16)}
\sum\limits_{i\in T}{{{x}_{i}}}\le D \mbox{\;\;\;(restrictions on the availability of weapons)},
\end{equation}
\begin{equation}\label{(17)}
{{x}_{i}}\in {{\mathbb{Z}}_{+}}, \;\;\;	i\in T,\ \mbox{(integer restrictions)},	
\end{equation}
where $T$ denotes the set of targets; ${{a}_{i}}$ - the cost of target $i\in T$; ${{p}_{i}}$ - the probability of hitting target $i$ with one weapon; ${{x}_{i}}$ - the number of weapons allocated to hit target $i$; $D$ is the available number of weapons; ${{\mathbb{Z}}_{+}}$ is the set of non-negative integers. Here 
\[{{\pi }_{i}}({{x}_{i}})={{\left( 1-{{p}_{i}} \right)}^{{{x}_{i}}}}={{\left( {{q}_{i}} \right)}^{{{x}_{i}}}}\] 
is the probability of not hitting target $i$ with the simultaneous independent use (volley) of ${{x}_{i}}$ weapons against target $i$; ${{q}_{i}}=1-{{p}_{i}}$. Problem \eqref{(15)}--\eqref{(17)} is also a nonlinear integer programming problem. For example, it can be solved by the dynamic programming method.

Another partial case of problem \eqref{(11)}--\eqref{(14)} is the problem of distributing single shots to targets, and one type of weapon can be used for different targets. Then the formulation \eqref{(11)}--\eqref{(14)} is supplemented with additional constraints:
\[
{{x}_{ij}}\in \left\{ 0,1 \right\},\;\;\;  i\in T,\ j\in W;\;\;\; 	\sum\nolimits_{j\in W}{{{x}_{ij}}}=1,\;\;\;	i\in T.
\] 
In addition, for ${{x}_{ij}}\in \left\{ 0,1 \right\}$, $\sum\nolimits_{j\in W}{{{x}_{ij}}}=1$, the equalities hold:
\[{{\left( 1-{{p}_{ij}} \right)}^{{{x}_{ij}}}}=1-{{x}_{ij}}{{p}_{ij}},\]     
\[{{\prod\nolimits_{j\in W}{\left( 1-{{p}_{ij}} \right)}}^{{{x}_{ij}}}}=\sum\nolimits_{j\in W}{\left( 1-{{p}_{ij}} \right){{x}_{ij}}}.\]   
Then problem \eqref{(11)}--\eqref{(14)} becomes a linear Boolean programming problem:
\begin{align*}
&\left[ \sum\nolimits_{i\in T,j\in W}{{{a}_{i}}}(1-{{p}_{ij}}){{x}_{ij}}=\sum\nolimits_{i\in T}{{{a}_{i}}\sum\nolimits_{j\in W}{\left( 1-{{p}_{ij}} \right){{x}_{ij}}}} \right]\to {{\min }_{\left\{ {{x}_{ij}} \right\}}},\\
&\sum\nolimits_{i\in T}{{{x}_{ij}}}\le {{D}_{j}},\;\;\;	j\in W;\\
&\sum\nolimits_{j\in W}{{{x}_{ij}}}=1,\;\;\;	i\in T;\\
&{{x}_{ij}}\in \left\{ 0,1 \right\},\;\;\;	i\in T,\ j\in W.
\end{align*}

\section{Method of approximate solution of problem \eqref{(11)}--\eqref{(14)}}

Let us return to problem \eqref{(11)}--\eqref{(14)}. Let us introduce additional variables ${{y}_{i}}=-\sum\nolimits_{j\in W}{{{x}_{ij}}\ln {{q}_{ij}}}$, $i\in T$. Then problem \eqref{(11)}--\eqref{(14)} can be equivalently reformulated as follows:
\begin{align}
&\sum\nolimits_{i\in T}{{{a}_{i}}}\,{{e}^{-{{y}_{i}}}}\to {{\min }_{\left\{ {{x}_{ij}},\,{{y}_{i}} \right\}}},\label{(18)}\\
&{{y}_{i}}+\sum\nolimits_{j\in W}{{{x}_{ij}}\ln \left( 1-{{p}_{ij}} \right)}=0,\;\;\;	i\in T, \label{(19)}\\
&\sum\limits_{i\in T}{{{x}_{ij}}}\le {{D}_{j}},\;\;\;	j\in W,\label{(20)}\\
&{{y}_{i}}\ge \ln \left( \varepsilon _{i}^{-1} \right),\;\;\;	i\in T,\label{(21)}\\
&{{x}_{ij}}\in {{\mathbb{Z}}_{+}}=\left\{ 0,1,2,... \right\},\;\;\;	i\in T,\ j\in W. \label{(22)}
\end{align}

Obviously, for the components 
\[{{y}_{i}}=\sum\nolimits_{j\in W}{{{b}_{ij}}{{x}_{ij}}},\;\;\;{{b}_{ij}}=-\ln \left( 1-{{p}_{ij}} \right),\] 
the optimal solution $\left\{ {{x}_{ij}},{{y}_{i}} \right\}$ from \eqref{(19)}, \eqref{(20)} the inequalities follow: 
\[{y}_{i}\le {{Y}_{i}}, \;\;\;\sum\nolimits_{i}{{{y}_{i}}}\le Y,\] 
where
\begin{align}
&{{y}_{i}}\le{{\max }_{{{x}_{ij}}\ge 0}}\left\{ \sum\nolimits_{j\in W}{{{b}_{ij}}{{x}_{ij}}}:\,\,\sum\limits_{i\in T}{{{x}_{ij}}}\le {{D}_{j}},\,\,j\in W \right\}\le \sum\nolimits_{j\in W}{{{b}_{ij}}{{D}_{j}}}=Y_i \label{(23)}\\
&\sum\nolimits_{i}{{{y}_{i}}}\le \sum\nolimits_{i}{\sum\nolimits_{j}{{{b}_{ij}}{{x}_{ij}}}}=\sum\nolimits_{j}{\sum\nolimits_{i}{{{b}_{ij}}{{x}_{ij}}}}\,\le \sum\nolimits_{j}{\sum\nolimits_{i}{\left( {{\max }_{i}}{{b}_{ij}} \right){{x}_{ij}}}}\nonumber\\
&\;\;\;\;\;\;\;\;\;\;=\sum\nolimits_{j}{\left( {{\max }_{i}}{{b}_{ij}} \right)\sum\nolimits_{i}{{{x}_{ij}}}}\,\le \sum\nolimits_{j}{\left( {{\max }_{i}}{{b}_{ij}} \right){{D}_{j}}}=Y. \label{(24)}
\end{align}
Then problem \eqref{(18)}--\eqref{(22)} can be reformulated as follows:
\begin{align}
&\sum\nolimits_{i\in T}{{{a}_{i}}}\,{{e}^{-{{y}_{i}}}}\to {{\min }_{\left\{ {{x}_{ij}},\,{{y}_{i}} \right\}}},
\label{(25)}\\
& {{y}_{i}}+\sum\nolimits_{j\in W}{{{x}_{ij}}\ln \left( 1-{{p}_{ij}} \right)}=0,\;\;\;	i\in T;\label{(26)}\\
& {{y}_{i}}+\sum\nolimits_{j\in W}{{{x}_{ij}}\ln \left( 1-{{p}_{ij}} \right)}=0,	i\in T;\label{(27)}\\
& \sum\nolimits_{i\in T}{{{y}_{i}}}\le Y;\label{(28)}\\
& 0\le \ln \left( \varepsilon _{i}^{-1} \right)={\underline{y_i}}\le {{y}_{i}}\le {{Y}_{i}},
\;\;\;	i\in T;	\label{(29)}\\
& {{x}_{ij}}\in {{\mathbb{Z}}_{+}},\;\;\;	i\in T,\ j\in W.\label{(30)}
\end{align}

Note that the objective function \eqref{(25)} is convex downwards, each exponent ${{e}^{-{{y}_{i}}}}$ in \eqref{(25)} can be approximated from below and above with arbitrary accuracy by convex piecewise linear maximum functions. Indeed, we partition each interval $\left[ {\underline{y_i}},{{Y}_{i}} \right]$ by points $\left\{ {{y}_{in}},n\in {{N}_{i}} \right\}$ such that ${\underline{y_i}}={{y}_{i1}}<{{y}_{i2}}<...<{{y}_{i\left| {{N}_{i}} \right|}}={{Y}_{i}}$. Then
\[\underset{n\in {{N}_{i}}}{\mathop{\max }}\,\left( {{e}^{-{{y}_{i\,n}}}}\left( 1+{{y}_{in}}-{{y}_{i}} \right) \right)\le {{e}^{-{{y}_{i}}}}\le \underset{n\in {{N}_{i}}\backslash {{Y}_{i}}}{\mathop{\max }}\,\left( {{e}^{-{{y}_{in}}}}+\frac{{{e}^{-{{y}_{i\,(n+1)}}}}-{{e}^{-{{y}_{i\,n}}}}}{{{y}_{i(n+1)}}-{{y}_{in}}}\left( {{y}_{i}}-{{y}_{in}} \right) \right).\]
The objective function \eqref{(25)} is approximated as follows:
\[\sum\nolimits_{i\in T}{{{a}_{i}}\underset{n\in {{N}_{i}}}{\mathop{\max }}\,\left( {{e}^{-{{y}_{i\,n}}}}\left( 1+{{y}_{in}}-{{y}_{i}} \right) \right)}\le \sum\nolimits_{i\in T}{{{a}_{i}}{{e}^{-{{y}_{i}}}}},\]
\[\sum\nolimits_{i\in T}{{{a}_{i}}{{e}^{-{{y}_{i}}}}}
\le \sum\nolimits_{i\in T}{{{a}_{i}}\underset{n\in {{N}_{i}}\backslash {{Y}_{i}}}{\mathop{\max }}\,\left( {{e}^{-{{y}_{in}}}}+\frac{{{e}^{-{{y}_{i\,(n+1)}}}}-{{e}^{-{{y}_{i\,n}}}}}{{{y}_{i(n+1)}}-{{y}_{in}}}\left( {{y}_{i}}-{{y}_{in}} \right) \right)}.\]
We introduce additional variables ${{z}_{i}}=\underset{n\in {{N}_{i}}}{\mathop{\max }}\,\left( {{e}^{-{{y}_{i\,n}}}}\left( 1+{{y}_{in}}-{{y}_{i}} \right) \right),\,\,\,\,\,i\in T,$ then problem \eqref{(25)}--\eqref{(30)} is approximated (from below) by the problem:
\[\sum\nolimits_{i\in T}{{{a}_{i}}}\,{{z}_{i}}\to {{\min }_{\left\{ {{x}_{ij}},\,{{y}_{i}},{{z}_{i}} \right\}}},\]
with constraints \eqref{(26)}--\eqref{(30)} and ${{z}_{i}}=\underset{n\in {{N}_{i}}}{\mathop{\max }}\,\left( {{e}^{-{{y}_{i\,n}}}}\left( 1+{{y}_{in}}-{{y}_{i}} \right) \right),\,\,i\in T.$ In turn, this problem is equivalent to the following linear partially integer optimization problem:
\begin{equation}\label{(31)}
\sum\nolimits_{i\in T}{{{a}_{i}}}\,{{z}_{i}}\to {{\min }_{\left\{ {{x}_{ij}},\,{{y}_{i}},{{z}_{i}} \right\}}},
\end{equation}
under constraints \eqref{(26)}--\eqref{(30)} and ${{z}_{i}}\ge \left( {{e}^{-{{y}_{i\,n}}}}\left( 1+{{y}_{in}}-{{y}_{i}} \right) \right),\,\,\,\,\,\,n\in {{N}_{i}},\,\,\,\,\,i\in T.$

Similarly, problem \eqref{(25)}--\eqref{(30)} is approximated (from above) by a linear partially integer optimization problem:
\[\sum\nolimits_{i\in T}{{{a}_{i}}}\,{{z}_{i}}\to {{\min }_{\left\{ {{x}_{ij}},\,{{y}_{i}},{{z}_{i}} \right\}}},\]
under constraints \eqref{(26)}--\eqref{(30)} and ${{z}_{i}}\ge {{e}^{-{{y}_{in}}}}+\frac{{{e}^{-{{y}_{i\,(n+1)}}}}-{{e}^{-{{y}_{i\,n}}}}}{{{y}_{i(n+1)}}-{{y}_{in}}}\left( {{y}_{i}-{y}_{in}} \right),\,\,\,\,\,\,n\in {{N}_{i}}\backslash {{Y}_{i}},\,\,\,\,\,i\in T.$

Another (stepwise) method of approximating the objective function \eqref{(25)} is as follows. Indeed, we divide each interval $\left[ \underline{y_i},{{Y}_{i}} \right]$ by points $\left\{ {{y}_{in}},n\in {{N}_{i}} \right\}$ such that ${\underline{y_i}}={{y}_{i1}}<{{y}_{i2}}<...<{{y}_{i\left| {{N}_{i}} \right|}}={{Y}_{i}}$. We introduce binary variables $\left\{ {{z}_{in}}\in \left\{ 0,1 \right\},n\in {{N}_{i}},\,i\in T \right\}$ and add the following conditions to the constraints \eqref{(26)}--\eqref{(30)}:

	${{z}_{in}}=\left\{ \begin{matrix}
   1, & {{y}_{i}}\in \left[ {{y}_{in}},{{y}_{i(n+1)}} \right),  \\
   0, & {{y}_{i}}\notin \left[ {{y}_{in}},{{y}_{i(n+1)}} \right),  \\
\end{matrix} \right.$ 	$n\in {{N}_{i}}\backslash {{Y}_{i}}$.

\noindent
Then the objective function \eqref{(25)} is approximated (from below) by the expression:
\[\sum\limits_{i\in T}{{{a}_{i}}}\sum\nolimits_{n\in {{N}_{i}}\backslash {{Y}_{i}}}{{{e}^{-{{y}_{in}}}}{{z}_{in}}}.\]

\section{The task of maximizing the probability of hitting all targets}
Consider the problem of optimal distribution of weapons among targets to maximize the probability of hitting all (or a fraction of) targets simultaneously:
\begin{align}
&\prod\nolimits_{i\in T}{\left( 1-{{\prod\nolimits_{j\in W}{\left( 1-{{p}_{ij}} \right)}}^{{{x}_{ij}}}} \right)}\to {{\max }_{\left\{ {{x}_{ij}} \right\}}},\label{(32)}\\
&\mbox{\;\;\;(maximizing the probability of damage),}\nonumber\\
&\sum\nolimits_{i\in T}{{{x}_{ij}}}\le {{D}_{j}},\;\;\;	j\in W, \mbox{\;\;\;(restrictions on the availability of weapons)}\label{(33)}\\
&	\sum\nolimits_{i\in T}{{{x}_{ij}}}\le {{D}_{j}},\;\;\;	j\in W, \mbox{\;\;\;(integer restrictions).} \label{(34)}	
\end{align}

First, instead of maximizing the function \eqref{(32)}, we will maximize its logarithm, and second, we will introduce new variables ${{y}_{i}}=-\sum\nolimits_{j\in W}{{{x}_{ij}}\ln {{q}_{ij}}}$. Then, problem \eqref{(32)}--\eqref{(34)}, taking into account constraints \eqref{(10)} and estimates \eqref{(23)}, \eqref{(24)}, takes the form:
\begin{align}
&  \sum\nolimits_{i\in T}{\ln \left( 1-{{e}^{-{{y}_{i}}}} \right)}\,\to {{\max }_{\left\{ {{x}_{ij}},\,{{y}_{i}} \right\}}};\label{(35)}\\
&  {{y}_{i}}+\sum\nolimits_{j\in W}{{{x}_{ij}}\ln \left( 1-{{p}_{ij}} \right)}=0,\;\;\;	i\in T,	\label{(36)}\\
&  \sum\limits_{i\in T}{{{x}_{ij}}}\le {{D}_{j}},	\;\;\;j\in W;\label{(37)}\\
&  \sum\nolimits_{i\in T}{{{y}_{i}}}\le Y;\label{(38)}\\
&  0\le \ln \left( \varepsilon _{i}^{-1} \right)={\underline{y_i}}\le {{y}_{i}}\le {{Y}_{i}},	
\;\;\;i\in T;\label{(39)}\\
& {{x}_{ij}}\in {{\mathbb{Z}}_{+}},\;\;\;	i\in T,\ j\in W.	 \label{(40)}
\end{align}

Note that the objective function \eqref{(35)} is concave (convex upwards). We will apply the approximation of function \eqref{(35)} from below and above by minimum functions.

Divide each interval $\left[ \underline{y_i},{{Y}_{i}} \right]$ by points $\left\{ {{y}_{in}},n\in {{N}_{i}} \right\}$ such that $\underline{y_i}={{y}_{i1}}<{{y}_{i2}}<...<{{y}_{i\left| {{K}_{i}} \right|}}={{Y}_{i}}$. Then,
\[\underset{n\in {{N}_{i}}}{\mathop{\min }}\,\left( \ln \left( 1-{{e}^{-{{y}_{in}}}} \right)+\frac{{{e}^{-{{y}_{i\,n}}}}}{1-{{e}^{-{{y}_{i}}}}}\left( {{y}_{i}}-{{y}_{in}} \right) \right)\ge \ln \left( 1-{{e}^{-{{y}_{i}}}} \right),\]
\[\ln \left( 1-{{e}^{-{{y}_{i}}}} \right)\ge \underset{n\in {{N}_{i}}\backslash {{Y}_{i}}}{\mathop{\min }}\,\left( \ln \left( 1-{{e}^{-{{y}_{in}}}} \right)+\frac{{{y}_{i}}-{{y}_{in}}}{{{y}_{i(n+1)}}-{{y}_{in}}}\ln \frac{1-{{e}^{-{{y}_{i\,(n+1)}}}}}{1-{{e}^{-{{y}_{i\,n}}}}} \right).\]
Thus, problem \eqref{(35)}--\eqref{(40)} is approximated (from above) by the following linear partial integer programming problem:
\[\sum\nolimits_{i\in T}{{{z}_{i}}}\,\to {{\max }_{\left\{ {{x}_{ij}},\,{{y}_{i}},{{z}_{i}} \right\}}}\]
with constraints \eqref{(36)}--\eqref{(40)} and additional constraints-inequalities:
\[{{z}_{i}}\le \ln \left( 1-{{e}^{-{{y}_{in}}}} \right)+\frac{{{e}^{-{{y}_{i\,n}}}}}{1-{{e}^{-{{y}_{i}}}}}\left( {{y}_{i}}-{{y}_{in}} \right),\;\;\;\,n\in {{N}_{i}},\,\,\,\,\,i\in T.\]
Problem \eqref{(35)}--\eqref{(40)} is approximated (from below) by the problem:
\[\sum\nolimits_{i\in T}{{{z}_{i}}}\,\to {{\max }_{\left\{ {{x}_{ij}},\,{{y}_{i}},{{z}_{i}} \right\}}}\]
with constraints \eqref{(36)}--\eqref{(40)} and additional constraints-inequalities:
\[{{z}_{i}}\le \ln \left( 1-{{e}^{-{{y}_{in}}}} \right)+\frac{{{y}_{i}}-{{y}_{in}}}{{{y}_{i(n+1)}}-{{y}_{in}}}\ln \frac{1-{{e}^{-{{y}_{i\,(n+1)}}}}}{1-{{e}^{-{{y}_{i\,n}}}}},
\;\;\;n\in {{N}_{i}}\backslash {{Y}_{i}},\,\,\,\,\,i\in T.\]

\section{Spatial problem of optimal allocation of weapons among targets}

The spatial problem of optimal distribution of weapons among targets (Spatial Weapon-Target Allocation, or SWTA problem) with additional indices $k\in K,\,\,l\in L$ is formulated as follows:
\begin{align*}
&\sum\nolimits_{k\in K,\,i\in T}{{{a}_{ki}}}\left( 1-{{\prod\nolimits_{l\in L,\,j\in W}{\left( 1-{{p}_{ki\,lj}} \right)}}^{{{x}_{ki\,lj}}}} \right)\to {{\max }_{\left\{ {{x}_{ki\,lj}} \right\}}},\\
&\sum\nolimits_{k\in K,\,i\in T}{{{x}_{ki\,lj}}}\le {{D}_{lj}},		\;\;\;l\in L,\,\,j\in W;\\
&{{x}_{ki\,lj}}\in {{\mathbb{Z}}_{+}},\;\;\;	k\in K,\,\,l\in L,\,\,i\in T,\ j\in W.
\end{align*}
where $K$ denotes the set of target positions, $L$ denotes the set of weapon positions, $T$ denotes the set of target types, among which there may be targets of the same type; $W$ is the set of weapon types; ${{a}_{ki}}$ is the cost of the target at position $k\in K$; ${{p}_{ki\,lj}}$ is the probability of hitting the target at position $k$ by a weapon from position $l$; ${{x}_{ki\,lj}}$ is the number of weapons of the type from position $l$ that are simultaneously (in a volley) used against the target at position $k$, ${{D}_{lj}}$ is the available number of weapons of the type at position $l$, is the set of non-negative integers. Here 
\[{{\pi }_{ki}}=\prod\nolimits_{l\in L,\,j\in W}{{{\left( 1-{{p}_{ki\,lj}} \right)}^{{{x}_{ki\,lj}}}}}=\prod\nolimits_{l\in L,\,j\in W}{{\left( {{q}_{ki\,lj}} \right)}^{{{x}_{ki\,lj}}}}\] 
is the probability of not hitting the target at position $k$ with the simultaneous independent use (volley) of ${{x}_{ki\,lj}}$, $l\in L,\,\,j\in W,$ of all types of weapons from all positions against the target at position $k$; ${{q}_{ki\,lj}}=1-{{p}_{ki\,lj}}$. If for some $k,\,i$ all ${{x}_{ki\,lj}}=0$, $l\in L,\,\,j\in W,$ then this goal $\left( k,i \right)$ is not attacked. The problem is a four-index nonlinear integer programming problem.

The average cost of hit targets alone may not be sufficient to assess the effectiveness of the weapons distribution, since some targets must be hit with almost certainty. This requirement is taken into account by additional constraints for the problem, that the probability of hitting a target $\left( k,i \right)$ must be at least $1-{{\varepsilon }_{ki}}$:
\[1-{{\prod\nolimits_{l\in L,j\in W}{\left( 1-{{p}_{ki\,lj}} \right)}}^{{{x}_{ki\,lj}}}}\ge 1-{{\varepsilon }_{ki}},\]
or
\[{{\prod\nolimits_{l\in L,j\in W}{\left( 1-{{p}_{ki\,lj}} \right)}}^{{{x}_{ki\,lj}}}}\le {{\varepsilon }_{ki}},\]
${{\varepsilon }_{ki}}\in \left( 0,1 \right]$,   $k\in K,\,\,i\in T.$

If ${{\varepsilon }_{ki}}=1$, then the corresponding constraint is satisfied automatically and may not be taken into account in the problem statement. These inequalities can be equivalently rewritten in linear form:
\[
\sum\nolimits_{l\in L,j\in W}{\left( \ln \left( 1-{{p}_{ki\,lj}} \right) \right)}\,{{x}_{ki\,lj}}\le \ln {{\varepsilon }_{ki}},\;\;\;k\in K,\,\,i\in T.
\]

Using the additional variables ${{y}_{ki}}=-\sum\nolimits_{l\in L,\,j\in W}{{{x}_{ki\,lj}}\ln {{q}_{ki\,lj}}}$, $k\in K,\,\,i\in T$, the problem can be reformulated as follows:
\begin{align*}
&  \sum\nolimits_{k\in K,i\in T}{{{a}_{ki}}}\,{{e}^{-{{y}_{ki}}}}\to {{\min }_{\left\{ {{x}_{ki\,lj}},\,{{y}_{ki}} \right\}}},\\
&  {{y}_{ki}}+\sum\nolimits_{l\in L,j\in W}{{{x}_{ki\,lj}}\ln \left( 1-{{p}_{ki\,lj}} \right)}=0,\;\;\;
k\in K,\,\,i\in T;\\
&  \sum\nolimits_{k\in K,\,i\in T}{{{x}_{ki\,lj}}}\le {{D}_{lj}},\;\;\;	l\in L,\,\,j\in W;\\
&  {{y}_{ki}}\ge \ln \left( \varepsilon _{ki}^{-1} \right),\;\;\;	k\in K,\,\,i\in T;\\
&  {{x}_{ki\,lj}}\in {{\mathbb{Z}}_{+}},\;\;\;	k\in K,\,\,l\in L,\,\,i\in T,\ j\in W.
\end{align*}

Obviously, for the components ${{y}_{ki}}=\sum\nolimits_{l\in L,j\in W}{{{b}_{ki\,lj}}{{x}_{ki\,lj}}}$, where ${{b}_{ki\,lj}}=\ln {{\left( {{q}_{ki\,lj}} \right)}^{-1}}$, the optimal solution $\left\{ {{x}_{ki\,lj}},{{y}_{ki}} \right\}$ satisfies the inequalities: $\ln \left( \varepsilon _{ki}^{-1} \right)\le {{y}_{ki}}\le {{Y}_{ki}}$, \[\sum\nolimits_{k\in K,i\in T}{{{y}_{ki}}}\le Y,\] 
where
\begin{eqnarray*}
{{y}_{ki}}&=&{{\max }_{{{x}_{ij}}\ge 0}}\left\{ \sum\nolimits_{l\in L,j\in W}{{{b}_{ki\,lj}}{{x}_{ki\,lj}}}:\,\,\sum\nolimits_{k\in K,\,i\in T}{{{x}_{ki\,lj}}}\le {{D}_{lj}},\,\,l\in L,\,\,j\in W \right\}\\
&\le& \sum\nolimits_{l\in L,j\in W}{{{b}_{ki\,lj}}{{D}_{lj}}}={{Y}_{ki}};
\end{eqnarray*}
\begin{eqnarray*}
 \sum\nolimits_{k\in K,i\in T}{{{y}_{ki}}}&\le& \sum\nolimits_{k,i}{\sum\nolimits_{l,j}{{{b}_{ki\,lj}}{{x}_{ki\,lj}}}}=\sum\nolimits_{l\in L,j\in W}{\sum\nolimits_{k\in K,i\in T}{{{b}_{ki\,lj}}{{x}_{ki\,lj}}}}\,\\
&\le& \sum\nolimits_{j}{\sum\nolimits_{i}{\left( {{\max }_{k,i}}{{b}_{ki\,lj}} \right){{x}_{ki\,lj}}}} \\ 
& =&\sum\nolimits_{l\in L,j\in W}{\left( {{\max }_{k,i}}{{b}_{ki\,lj}} \right)\sum\nolimits_{k\in K,i\in T}{{{x}_{ki\,lj}}}}\,\\
&\le& \sum\nolimits_{l\in L,j\in W}{\left( \underset{k\in K,i\in T}{\mathop{\max }}\,{{b}_{ki\,lj}} \right){{D}_{lj}}}=Y. 
\end{eqnarray*}
Then the spatial problem can be reformulated as follows:
\begin{align*}
&  \sum\nolimits_{k\in K,i\in T}{{{a}_{ki}}}\,{{e}^{-{{y}_{ki}}}}\to {{\min }_{\left\{ {{x}_{ki\,lj}},\,{{y}_{ki}} \right\}}},\\
&  {{y}_{ki}}+\sum\nolimits_{l\in L,j\in W}{{{x}_{ki\,lj}}\ln \left( 1-{{p}_{ki\,lj}} \right)}=0,	\;\;\;
k\in K,\,\,i\in T;\\
&  \sum\nolimits_{k\in K,\,i\in T}{{{x}_{ki\,lj}}}\le {{D}_{lj}},	\;\;\;l\in L,\,\,j\in W;\\
&  \sum\nolimits_{k\in K,i\in T}{{{y}_{ki}}}\le Y;\\
&  0\le \ln \left( \varepsilon _{ki}^{-1} \right)=\underline{y_{ki}}\le {{y}_{ki}}\le {{Y}_{ki}},	\;\;\;k\in K,\,\,i\in T;\\
&  	{{x}_{ki\,lj}}\in {{\mathbb{Z}}_{+}},\;\;\;k\in K,\,\,l\in L,\,\,i\in T,\ j\in W.
\end{align*}

To solve this problem, we can apply the same approximation approach as in the previous section.

\section{Critical Infrastructure Protection Modeling}

Modeling the protection of critical infrastructure facilities from air attacks consists in the optimal distribution of defensive weapons among the facilities. In this case, the objective function can take into account not only the cost of destroyed targets, but also the lost cost of the facility and, possibly, the cost without the weapons used back. Such models (in a game setting) are considered in \cite{norkin2025simulation,norkin2025allocation}. The task of the defending party, with a known attack plan, is to minimize the damage from the attack by optimally distributing defensive resources.

{\bf Model 1} (several homogeneous units attack several objects protected by homogeneous means of defense):
\begin{equation}\label{(41)}
U(a,x)=\sum\nolimits_{k\in K}{{{c}_{k}}\left( 1-{{\left( 1-{{\left( 1-p \right)}^{{{x}_{k}}}} \right)}^{{{a}_{k}}}} \right)}
\end{equation}
(Expected losses to be minimized by $\left\{ {{x}_{k}} \right\}$),
\begin{equation}\label{(42)}
\sum\limits_{i=1}^{K}{{{x}_{i}}}\le W \mbox{\;\;\;protective resources (units)},
\end{equation}
where
$K$ is set of defense objects;

${{c}_{k}}$ is value of object $k$;

$p$ is probability of defeating an attacking unit by one defense unit;

${{a}_{k}}$ is given (known) number of attacking units at object $k$;

${{x}_{k}}$ is number of weapons for defense at object $k$;

${{q}_{k}}=1-{{\left( 1-p \right)}^{{{x}_{k}}}}$ is probability of defeating an attacking unit ${{x}_{k}}$ by defense units; 

$p({{a}_{k}},{{x}_{k}})=\left( 1-{{\left( {{q}_{k}} \right)}^{{{x}_{k}}}} \right)=1-{{\left( 1-{{\left( 1-p \right)}^{{{x}_{k}}}} \right)}^{{{a}_{k}}}}$ is the probability that at least one unit that makes the attack will hit object $k$;

$a=\left( {{a}_{1}},...,{{a}_{K}} \right)$; $x=\left( {{x}_{1}},...,{{x}_{K}} \right)$.

In the case where all ${{a}_{k}}=1$, the objective function is simplified: $U(a,x)=\sum\nolimits_{k\in K}{{{c}_{k}}{{\left( 1-p \right)}^{{{x}_{k}}}}}$. To take into account the cost of spent defense resources, the objective function takes the form 
$U(a,x)=\sum\nolimits_{k\in K}{\left( {{c}_{k}}\left( 1-{{\left( 1-{{\left( 1-p \right)}}^{{{x}_{k}}} \right)}^{{{a}_{k}}}} \right)+b{{x}_{k}} \right)}$, 
where $b$ is the cost of one lost weapon unit. These problems can be solved by the dynamic programming method, since the objective functions are separable.

If the relative value of objects is difficult to estimate, then the problem can be reformulated as the problem of maximizing the damage of all attack units:
\begin{align}
&P(a,x)=\prod\nolimits_{k\in K}\left( 1-{{\left( 1-p \right)}^{{{x}_{k}}}} \right)^{{{a}_{k}}}
\to {\max }_{\left\{ {x}_{i} \right\}},\label{(43)}\\
&\sum\nolimits_{k\in K}{{x}_{k}}\le \left| W \right|,\label{(44)}\\
&{{x}_{k}}\in \mathbb{Z}_{+},	\;\;\;k\in K.\label{(45)}
\end{align}

Problem \eqref{(43)}--\eqref{(45)} can be solved by dynamic programming, and can also be reduced to a problem of partial integer programming. We prologarithmize the objective function in \eqref{(43)} and introduce new variables, then we obtain an equivalent problem:
\begin{align}
&  \ln P(a,x)=\sum\nolimits_{k\in K}{a}_{k}\ln \left( 1-{e}^{-{{y}_{k}}} \right)\to {\max }_{\left\{ {x}_{k},{y}_{k} \right\}},\label{(46)}\\
&  {y}_{k}+{x}_{k}\ln (1-p)=0, 	\;\;\;{x}_{k}\in \mathbb{Z}_{+},	\;\;\;k\in K.	\label{(47)}
\end{align}

The functions $\ln \left( 1-{{e}^{-{{y}_{k}}}} \right)$ are convex upwards and therefore can be approximated by piecewise linear convex upwards functions. Therefore, the problem can also be reduced to a linear partial integer programming problem.

{\bf Model 2} (combined attack on one object). In this model, it is assumed that the object is attacked by $a=\left\{ {{a}_{j}},\,\,j\in T \right\}$ attack units, and protected by ${{D}_{i}},\,\,i\in W,$ defense units, ${{p}_{ij}}$ is the probability of defeating one attack unit by one defense unit. The problem of maximizing the probability of defeating all attack units by optimizing the defense plan $x=\left\{ {{x}_{ij}},\,\,i\in W,\,\,j\in T \right\}$ has the form:
\begin{align*}
&P(a,x)=\prod\limits_{j\in T}{{{\left( 1-\prod\nolimits_{i\in W}{{{\left( 1-{{p}_{ij}} \right)}^{{{x}_{ij}}}}} \right)}^{{{a}_{j}}}}}\to {{\max }_{\left\{ {{x}_{ij}} \right\}}},\\
&\sum\limits_{j\in T}{{{x}_{ij}}}\le {{D}_{i}},\,\,\,\,\,\,\,i\in W\,;\\
& {{x}_{ij}}\in {{\mathbb{Z}}_{+}},\;\;\;	i\in T,\ j\in W.	
\end{align*}

To solve this problem, we can also apply the reduction method to a convex partially integer optimization problem. To do this, first we logarithmize the objective function: $\ln P(a,x)=\sum\limits_{j\in T}{{{a}_{j}}\left( 1-\prod\nolimits_{i\in W}{{{\left( 1-{{p}_{ij}} \right)}^{{{x}_{ij}}}}} \right)}$, and then using the new variables ${{y}_{i}}=-\sum\nolimits_{j\in W}{{{x}_{ij}}\ln \left( 1-{{p}_{ij}} \right)}$, $i=1,...,T$, we obtain the equivalent problem:
\begin{align*}
&   \ln P(a,x)=\sum\limits_{j\in T}{{{a}_{j}}\left( 1-{{e}^{-{{y}_{j}}}} \right)}\to {{\max }_{\left\{ {{x}_{ij}},{{y}_{i}} \right\}}},\\
&   {{y}_{i}}+\sum\nolimits_{j\in W}{{{x}_{ij}}\ln \left( 1-{{p}_{ij}} \right)}=0,\;\;\; i\in T;\\
&   \sum\limits_{j\in T}{{{x}_{ij}}}\le {{D}_{i}},\,\,\,\,\,\,\,i\in W\,; 	\\
&   {{x}_{ij}}\in {{\mathbb{Z}}_{+}},\;\;\;	i\in T,\ j\in W.
\end{align*}

This problem with additional constraints \eqref{(23)}--\eqref{(24)} has an upward convex objective function and is similar to problem \eqref{(35)}--\eqref{(40)}, therefore it can be solved by the reduction method to linear partial integer programming.

{\bf Model 3} (combined protection of several objects from a combined attack \cite{norkin2025simulation}. Suppose that the attack plan $\left\{ {{a}_{kj}},\,\,k\in K,\,\,j\in T \right\}$ on objects $k\in K$ by attack means $j\in T$ is known. Then the optimization problem of the protection plan $\left\{ {{x}_{ki}} \right\}$ has the form:

$U(a,x)=\sum\limits_{k\in K}{{{c}_{k}}\left( 1-\prod\nolimits_{j\in T}{{{\left( 1-\prod\nolimits_{i\in W}{{{\left( 1-{{p}_{ij}} \right)}^{{{x}_{ki}}}}} \right)}^{{{a}_{kj}}}}} \right)}$ (Expected losses to be minimized by $\left\{ {{x}_{ki}} \right\}$),

$\sum\limits_{k\in K}{{{x}_{ki}}}\le {{D}_{i}},\,\,\,\,\,\,\,i\in W,$ are defense resources,

$K$ is the set of objects to protect, $k\in K$,

${{c}_{k}}$ is the  value of object $k$,

$W$ is the set of types of defensive units,

$T$ is the set of types of attacking units,

${{p}_{ij}}$ is the probability of defeating an attacking unit of type $j$ by one defensive unit of type $i$,

${{a}_{kj}}$ is the number of attacking units of type $j$ per object $k$,

${{x}_{ki}}\in {{\mathbb{Z}}_{+}}$ is the number of protective units of type $i$ on object $k$,

${{a}_{k}}=\left( {{a}_{k1}},...,{{a}_{kT}} \right)$,

 ${{x}_{k}}=\left( {{x}_{11}},...,{{x}_{kW}} \right)$,

${{q}_{kj}}=1-\prod\nolimits_{i\in W}{{{\left( 1-{{p}_{ij}} \right)}^{{{x}_{ki}}}}}$ is the  probability of defeating one attacking unit of type $j$ ${{x}_{ki}}$ by means of defense units at object $k$,

$p({{a}_{k}},{{x}_{k}})=\left( 1-\prod\nolimits_{j\in T}{{{q}_{kj}}} \right)=1-\prod\nolimits_{j\in T}{{{\left( 1-\prod\nolimits_{i\in W}{{{\left( 1-{{p}_{ij}} \right)}^{{{x}_{ki}}}}} \right)}^{{{a}_{kj}}}}}$ is the probability that at least one attacking unit will hit object $k$.

If the goal of protection is to maximize the probability of defeating all attack units, then the objective function takes the form:
$$P(a,x)=\prod\limits_{k\in K}{\prod\nolimits_{j\in T}{{{\left( 1-\prod\nolimits_{i\in W}{{{\left( 1-{{p}_{ij}} \right)}^{{{x}_{ki}}}}} \right)}^{{{a}_{kj}}}}}}\to {{\max }_{\left\{ {{x}_{ki}} \right\}}}.$$

\section{Accounting the uncertainty regarding the attack plan}

Suppose that for model 3 the attack plan $\left\{ {{a}_{kj}},\,\,k\in K,\,\,j\in T \right\}$ on objects $k\in K$ by means $j\in T$ is not known for certain, but $s\in S$ attack scenarios ${{a}_{s}}=\left\{ {{a}_{kjs}},\,\,k\in K,\,\,j\in T \right\}$ with corresponding probabilities ${{q}_{s}} are known.$ Then the optimization problem of the protection plan $\left\{ {{x}_{ki}} \right\}$ has the form:

$U(a,x)=\sum\limits_{s\in S}{{{q}_{s}}\sum\limits_{k\in K}{{{c}_{k}}\left( 1-\prod\nolimits_{j\in T}{{{\left( 1-\prod\nolimits_{i\in W}{{{\left( 1-{{p}_{ij}} \right)}^{{{x}_{ki}}}}} \right)}^{{{a}_{kjs}}}}} \right)}}$ (Expected losses to be minimized by $\left\{ {{x}_{ki}} \right\}$), 

$\sum\limits_{k\in K}{{{x}_{ki}}}\le {{D}_{i}},\,\,\,\,\,\,\,i\in W,$ (Defense resources).

\section{Numerical example}

Let us consider a numerical example of solving the problem \eqref{(11)} -- \eqref{(14)} by reducing it to a linear problem of partial integer programming \eqref{(31)}. The initial problem \eqref{(11)} -- \eqref{(14)} has the following parameters: the number of targets $\left| T \right|=14$, the cost of the targets, $a=\left( 0.1,0.1,0.1,0.1,0.1,0.1,0.1,0.1,0.1,0.1,0.5,1.0,1.5,1.5,1.5 \right)$, the number of types of weapons $\left| W \right|=3$, number of weapons $\left\{ {{D}_{j}},\,j=1,2,3 \right\}=\left\{ 5,30,30 \right\}$, probability of hitting targets 
\begin{center}
$\varepsilon =\left( 0.3,0.3,0.3,0.3,0.3,0.3,0.3,0.3,0.3,0.3,0.4,0.3,0.8,0.8,0.8 \right)$, 
\end{center}
(transposed) weapon efficiency matrix ${P}'=\left\{ {{p}_{ji}} \right\}$:

\begin{small}                
\[\left( \begin{array}{@{}cccccccccccccc@{}}
   0.7 & 0.7 & 0.7 & 0.7 & 0.7 & 0.7 & 0.7 & 0.7 & 0.7 & 0.8 & 0.5 & 0.2 & 0.2 & 0.2  \\
   0.5 & 0.5 & 0.5 & 0.5 & 0.5 & 0.5 & 0.5 & 0.5 & 0.5 & 0 & 0 & 0 & 0 & 0  \\
   0.2 & 0.2 & 0.2 & 0.2 & 0.2 & 0.2 & 0.2 & 0.2 & 0.2 & 0.1 & 0.03 & 0.001 & 0.001 & 0.001  \\
\end{array} \right).
\]
\end{small}

\smallskip
The corresponding problem \eqref{(31)} with $Y=35.5359$, \\$\left\{ {{Y}_{i}},\,i\in T \right\}=\left\{ 33.5086,33.5086,...,33.5086,11.2080,4.3785,1.1457,1.1457,1.1457 \right\}$,\\
$\left\{ \left| {{N}_{i}} \right|,\,\,i=1,...,T \right\}=\left\{ 35,35,35,35,35,35,35,35,13,6,3,3,3,3 \right\}$ \\(approximation of functions ${{e}^{-{{y}_{i}}}}$ by tangents at integer points), had the following dimensions: the number of continuous variables $2\left| T \right|=28$ , the number of discrete variables $\left| T \right|\times \left| W \right|=14\times 3=42$ , the number of inequality constraints 347, the number of equality constraints 14, the number of two bound constraints on the variables 70. The solution to problem \eqref{(11)} -- \eqref{(14)}, obtained by the standard intlinprog program of the Matlab system, has the form:

\begin{small}
${x}'=\left\{ {{x}_{ji}} \right\}=\left( \begin{array}{@{}cccccccccccccc@{}}
   0 & 0 & 0 & 0 & 0 & 0 & 0 & 0 & 0 & 0 & 1 & 1 & 2 & 1  \\
   4 & 4 & 4 & 3 & 3 & 3 & 3 & 3 & 3 & 0 & 0 & 0 & 0 & 0  \\
   0 & 0 & 0 & 0 & 0 & 0 & 0 & 0 & 0 & 13 & 17 & 0 & 0 & 0  \\
\end{array} \right)$
\end{small}

\noindent
with the probability vector $\left( {{\text{ }\!\!\Pi\!\!\text{ }}_{j\in W}}{{\left( 1-{{p}_{ij}} \right)}^{{{x}_{ij}}}},\,\,i\in T \right)$ of not hitting the targets:
\begin{small}
$\left( 0.0625,\,0.0625,\,0.0625,\,0.125,\,0.125,\,0.125,\,0.125,\,0.125,\,0.125,\,0.2542,\,0.2979,\,0.8,\,0.64,\,0.8 \right)$.
\end{small}
The calculation time on an HP Pavilion Laptop 15 with an 11th Gen Intel(R) Core(TM) i5-1135G7 processor and 16.0 GB of RAM was $t\approx 0.185$sec.

\section{Conclusions}

The paper shows that known WTA models of optimal allocation of weapons by targets can be used for planning defense against air attacks. These optimization models are supplemented by constraints on the risk of not hitting targets (or losing the object of protection). The models belong to the class of nonlinear integer optimization problems. An approach is proposed for the approximate solution of these problems by reducing them to problems of linear partial integer programming of high dimension with control of the approximation accuracy. Both cost and probabilistic settings
can be solved by this approach.

\smallskip
{\bf Acknowledgment.}
The paper study was carried out within the project 2025.07/0403 ''Adaptive stochastic optimization methods for solving machine learning problems, statistical estimation, and control of dynamic systems under adversarial conditions and uncertainty'' (2026--2028) funded by the National Research Foundation of Ukraine (NRFU).

\bibliographystyle{unsrt}
\bibliography{Norkin_WTA_ref.bib}
\end{document}